\documentclass[10pt,reqno]{amsart}
\usepackage{amsmath,amsthm,amsfonts,amssymb,amscd}
\usepackage[latin1]{inputenc}
\usepackage{psfrag,pslatex}
\usepackage{epsfig}

\newcommand{\dist}{\operatorname{dist}}
\newcommand{\const}{\operatorname{const}}
\newcommand{\inter}{\operatorname{int}}

\newcommand{\Per}{\operatorname{Per}}
\newcommand{\ang}[1]{\operatorname{angle\,}(#1)}

\newcommand{\vep}{\varepsilon}

\def \CC {{\mathbb C}}
\def \DD {{\mathbb D}}

\def \NN {{\mathbb N}}

\def \RR {{\mathbb R}}
\def \SS {{\mathbb S}}
\def \TT {{\mathbb T}}

\def \ZZ {{\mathbb Z}}

\def \cF {{\mathcal F}}

\def \cP {{\mathcal P}}

\def \cR {{\mathcal R}}

\def \cV {{\mathcal V}}

\theoremstyle{remark}

\begin{document}

\title{Hyperbolic Dynamical Systems}

\author{Vitor Ara\'ujo}
\address{CMUP, Rua do Campo Alegre 687, 4169-007 Porto, Portugal}
\email{vdaraujo@fc.up.pt}
\address{IM-UFRJ, C.P. 68.530, CEP 21.945-970, Rio de Janeiro, Brazil}
\email{vitor.araujo@im.ufrj.br}

\author{Marcelo Viana}
\address{IMPA, Estrada Dona Castorina 110, CEP 22460-320 Rio de Janeiro, Brasil}
\email{viana@impa.br}

\date{\today}

\maketitle

\setcounter{tocdepth}{1}\tableofcontents

\setcounter{tocdepth}{1}\tableofcontents

\section*{Glossary}

\subsection*{Homeomorphism, diffeomorphism}

A \emph{homeomorphism} is a continuous map $f:M\to N$ which is
one-to-one and onto, and whose inverse $f^{-1}:N\to M$ is
also continuous. It may be seen as a global continuous change
of coordinates. We call $f$ a \emph{diffeomorphism} if, in addition,
both it and its inverse are smooth. When $M=N$, the iterated
$n$-fold composition $f\circ\stackrel{n}{\dots}\circ f$ is
denoted by $f^n$. By convention, $f^0$ is the identity map,
and $f^{-n}=(f^n)^{-1}=(f^{-1})^n$ for $n\ge0$.

\subsection*{Smooth flow}

A flow $f^t:M\to M$ is a family of diffeomorphisms depending in a
smooth fashion on a parameter $t\in\RR$ and satisfying
$f^{s+t}=f^s\circ f^t$ for all $s$, $t\in\RR$. This property
implies that $f^0$ is the identity map. Flows usually arise as
solutions of autonomous differential equations: let $t\mapsto
\phi^t(v)$ denote the solution of
\begin{equation}\label{eq.ode}
\dot X = F(X), \quad X(0)=v
\end{equation}
and assume solutions are defined for all times; then the family $\phi^t$
thus defined is a flow (at least as smooth as the vector field $F$ itself).
The vector field may be recovered from the flow, through the relation
$F(X) = \partial_t \phi^t(X)\mid_{t=0}$.

\subsection*{$C^k$ topology}

Two maps admitting continuous derivatives are said to be
\emph{$C^1$-close} if they are uniformly close, and so are their
derivatives. More generally, given any $k\ge 1$, we say that two
maps are \emph{$C^k$-close} if they admit continuous derivatives
up to order $k$, and their derivatives of order $i$ are uniformly
close, for every $i=0, 1, \ldots, k$. This defines a topology in
the space of maps of class $C^k$.

\subsection*{Foliation}

A foliation is a partition of a subset of the ambient space into smooth
submanifolds, that one calls leaves of the foliation, all with the same
dimension and varying continuously from one point to the other.
For instance, the trajectories of a vector field $F$, that is, the
solutions of equation \eqref{eq.ode}, form a $1$-dimensional foliation
(the leaves are curves) of the complement of the set of zeros of $F$.
The main examples of foliations in the context of this work are the
families of stable and unstable manifolds of hyperbolic sets.

\subsection*{Attractor}

A subset $\Lambda$ of the ambient space $M$ is
\emph{invariant} under a transformation $f$ if
$f^{-1}(\Lambda)=\Lambda$, that is, a point is in $\Lambda$
if and only if its image is. $\Lambda$ is invariant under a
flow if it is invariant under $f^t$ for all $t\in\RR$. An
\emph{attractor} is a compact invariant subset $\Lambda$
such that the trajectories of all points in a neighborhood
$U$ converge to $\Lambda$ as times goes to infinity, and
$\Lambda$ is \emph{dynamically indecomposable} (or
\emph{transitive}): there is some trajectory dense in
$\Lambda$. {Sometimes one asks convergence only for points
  in some ``large'' subset of a neighborhood $U$ of
  $\Lambda$, and dynamical indecomposability can also be
  defined in somewhat different ways. However, the
  formulations we just gave are fine in the uniformly
  hyperbolic context.}

\subsection*{Limit sets}

The \emph{$\omega$-limit set} of a trajectory $f^n(x)$, $n\in\ZZ$
is the set $\omega(x)$ of all accumulation points of the
trajectory as time $n$ goes to $+\infty$. The \emph{$\alpha$-limit
set} is defined analogously, with $n\to-\infty$. The corresponding
notions for continuous time systems (flows) are defined
analogously. The \emph{limit set} $L(f)$ (or $L(f^t)$, in the flow
case) is the closure of the union of all '$\omega$-limit and all
$\alpha$-limit sets. The \emph{non-wandering set} $\Omega(f)$ (or
$\Omega(f^t)$, in the flow case) is that set of points such that
every neighborhood $U$ contains some point whose orbit returns to
$U$ in future time (then some point returns to $U$ in past time as
well). When the ambient space is compact all these sets are
non-empty. Moreover,the limit set is contained in the
non-wandering set.

\subsection*{Invariant measure}

A probability measure $\mu$ in the ambient space $M$ is
\emph{invariant} under a transformation $f$ if
$\mu(f^{-1}(A))=\mu(A)$ for all measurable subsets $A$. This means
that the ``events'' $x\in A$ and $f(x)\in A$ have equally
probable. We say $\mu$ is invariant under a flow if it is
invariant under $f^t$ for all $t$.  An invariant probability
measure $\mu$ is \emph{ergodic} if every invariant set $A$ has
either zero or full measure. An equivalently condition is that
$\mu$ can not be decomposed as a convex combination of invariant
probability measures, that is, one can not have
$\mu=a\mu_1+(1-a)\mu_2$ with $0<a<1$ and $\mu_1$, $\mu_2$
invariant.

\section*{Definition}

In general terms, a smooth dynamical system is called hyperbolic
if the tangent space over the asymptotic part of the phase space
splits into two complementary directions, one which is contracted 
and the other which is expanded under the action of the system.
In the classical, so-called uniformly hyperbolic case, the 
asymptotic part of the phase space is embodied by the limit set
and, most crucially, one requires the expansion and contraction
rates to be uniform. Uniformly hyperbolic systems are now fairly
well understood. They may exhibit very complex behavior which,
nevertheless, admits a very precise description. Moreover, 
uniform hyperbolicity is the main ingredient for characterizing
structural stability of a dynamical system. Over the years the 
notion of hyperbolicity was broadened (non-uniform hyperbolicity)
and relaxed (partial hyperbolicity, dominated splitting) to 
encompass a much larger class of systems, and has become a 
paradigm for complex dynamcial evolution.

\section{Introduction}

The theory of uniformly hyperbolic dynamical systems was initiated
in the 1960's (though its roots stretch far back into the 19th
century) by S. Smale, his students and collaborators, in the west,
and D. Anosov, Ya. Sinai, V. Arnold, in the former Soviet Union.
It came to encompass a detailed description of a large class of
systems, often with very complex evolution. Moreover, it provided
a very precise characterization of structurally stable dynamics,
which was one of its original main goals.

The early developments were motivated by the problem of characterizing
structural stability of dynamical systems, a notion that had been
introduced in the 1930's by A. Andronov and L. Pontryagin.
Inspired by the pioneering work of M. Peixoto on circle maps and
surface flows, Smale introduced a class of \emph{gradient-like} 
systems, having a finite number of periodic orbits, which should be
structurally stable and, moreover, should constitute the majority
(an open and dense subset) of all dynamical systems.
Stability and openness were eventually established, in the thesis
of J. Palis. However, contemporary results of M. Levinson, based on
previous work by M. Cartwright and J. Littlewood, provided examples
of open subsets of dynamical systems all of which have an infinite 
number of periodic orbits.

In order to try and understand such phenomenon, Smale introduced a
simple geometric model, the now famous "horseshoe map",
for which infinitely many periodic orbits exist in a robust way.
Another important example of structurally stable system which is
not gradient like was R. Thom's so-called "cat map".  
The crucial common feature of these models is hyperbolicity:
the tangent space at each point splits into two complementar
directions such that the derivative contracts one of these
directions and expands the other, at uniform rates. 

In global terms, a dynamical system is called \emph{uniformly
hyperbolic}, or Axiom A, if its limit set has this hyperbolicity
property we have just described. The mathematical theory of such
systems, which is the main topic of this paper, is now well
developped and constitutes a main paradigm for the behavior of
"chaotic" systems. In our presentation we go from local aspects
(linear systems, local behavior, specific examples) to the global
theory (hyperbolic sets, stability, ergodic theory). In the final
sections we discuss several important extensions (strange
attractors, partial hyperbolicity, non-uniform hyperbolicity) that
have much broadened the scope of the theory.

\section{Linear systems}\label{sec:linear-systems}

Let us start by introducing the phenomenon of hyperbolicity
in the simplest possible setting, that of linear
transformations and linear flows. Most of what we are going
to say applies to both discrete time and continuous time
systems in a fairly analogous way, and so at each point we
refer to either one setting or the other. In depth
presentations can be found in e.g. \cite{PM82} and
\cite{KH95}.

The general solution of a system of linear ordinary
differential equations
$$
\dot X=AX,\quad X(0)=v
$$
where $A$ is a constant $n\times n$ real matrix and $v\in\RR^n$ is
fixed,  is given by
$$
X(t)=e^{tA}\cdot v, \quad t\in\RR,
$$
where $e^{tA}=\sum_{n=0}^\infty (tA)^n/n!$. The linear flow is
called  \emph{hyperbolic} if $A$ has no eigenvalues on the
imaginary axis. Then the \emph{exponential} matrix $e^{A}$ has no
eigenvalues with norm $1$. This property is very important for a
number of reasons.

\subsection*{Stable and unstable spaces}

For one thing it implies that all solutions have well-defined
asymptotic behavior: they either converge to zero or diverge to
infinity as time $t$ goes to $\pm\infty$. More precisely, let
\begin{itemize}
\item $E^s$ (\emph{stable subspace}) be the subspace of $\RR^n$ spanned
by the generalized eigenvector associated to eigenvalues of $A$
with negative real part.
\item $E^u$ (\emph{unstable subspace}) be the subspace of $\RR^n$ spanned
by the generalized eigenvector associated to eigenvalues of $A$
with positive real part
\end{itemize}
Then these subspaces are complementary, meaning that
$\RR^n=E^s\oplus E^u$, and every solution $e^{tA}\cdot v$ with
$v\not\in E^s\cup E^u$ diverges to infinity both in the future and
in the past. The solutions with $v\in E^s$ converge to zero as
$t\to+\infty$ and go to infinity as $t\to-\infty$, and analogously
when $v\in E^u$, reversing the direction of time.

\subsection*{Robustness and density}

Another crucial feature of hyperbolicity is \emph{robustness}: any
matrix that is close to a hyperbolic one, in the sense that
corresponding coefficients are close, is also hyperbolic. The
stable and unstable subspaces need not coincide, of course, but
the dimensions remain the same. In addition, hyperbolicity if
\emph{dense}: any matrix is close to a hyperbolic one. That is
because, up to arbitrarily small modifications of the
coefficients, one may force all eigenvalues to move out of the
imaginary axis.

\subsection*{Stability, index of a fixed point}

In addition to robustness, hyperbolicity also implies
\emph{stability}: if $B$ is close to a hyperbolic matrix $A$, in
the sense we have just described, then the solutions of $\dot
X=BX$ have essentially the same behavior as the solutions of $\dot
X=A X$. What we mean by ``essentially the same behavior'' is that
there exists a global continuous change of coordinates, that is, a
homeomorphism $h:\RR^n\to\RR^n$, that maps solutions of one system
to solutions of the other, preserving the time parametrization:
$$
h\big(e^{tA}\cdot v\big) = e^{t B}\cdot h(v)
 \quad\text{for all} \quad t\in\RR.
$$
More generally, two hyperbolic linear flows are conjugated by a
homeomorphism $h$ if and only if they have the same \emph{index},
that is, the same number of eigenvalues with negative real part.
In general, $h$ can not be taken to be a diffeomorphism: this is
possible if and only if the two matrices $A$ and $B$ are obtained
from one another via a change of basis. Notice that in this case
they must have the same eigenvalues, with the same multiplicities.

\subsection*{Hyperbolic linear flows}

There is a corresponding notion of hiperbolicity for discrete time
linear systems
$$
X_{n+1}=CX_n, \quad X_0=v
$$
with $C$ a $n\times n$ real matrix. Namely, we say the system is
\emph{hyperbolic} if $C$ has no eigenvalue in the unit circle.
Thus a matrix $A$ is hyperbolic in the sense of continuous time
systems if and only if its exponential $C=e^A$ is hyperbolic in
the sense of discrete time systems. The previous observations
(well-defined behavior, robustness, denseness and stability)
remain true in discrete time. Two hyperbolic matrices are
conjugate by a homeomorphism if and only if they have the same
index, that is, the same number of eigenvalues with norm less than
$1$, and they both either preserve or reverse orientation.

\section{Local theory}\label{sec:local-theory}

Now we move on to discuss the behavior of non-linear systems close
to fixed or, more generally, periodic trajectories. By non-linear
system we understand the iteration of a diffeomorphism $f$, or the
evolution of a smooth flow $f^t$, on some manifold $M$. The
general philosophy is that the behavior of the system close to a
hyperbolic fixed point very much resembles the dynamics of its
linear part.

A fixed point $p\in M$ of a diffeomorphism $f:M\to M$ is called
\emph{hyperbolic} if the linear part $Df_p:T_p M\to T_pM$ is a
hyperbolic linear map, that is, if $Df_p$ has no eigenvalue with
norm $1$. Similarly, an equilibrium point $p$ of a smooth vector
field $F$ is \emph{hyperbolic} if the derivative $DF(p)$ has no
pure imaginary eigenvalues.

\subsection*{Hartman-Grobman theorem}

This theorem asserts that if $p$ is a hyperbolic fixed point of
$f:M\to M$ then there are neighborhoods $U$ of $p$ in $M$ and $V$
of $0$ in the tangent space $T_p M$ such that we can find a
homeomorphism $h:U\to V$ such that
$$
h \circ f = Df_p \circ h
$$
whenever the composition is defined. This property means that $h$
maps orbits of $Df(p)$ close to zero to orbits of $f$ close to $p$.
We say that $h$ is a (local) \emph{conjugacy} between the non-linear
system $f$ and its linear part $Df_p$. There is a corresponding
similar theorem for flows near a hyperbolic equilibrium. In either
case, in general $h$ can not be taken to be a diffeomorphism.

\subsection*{Stable sets}

The \emph{stable set} of the hyperbolic fixed point $p$ is defined by
\begin{align*}
  W^s(p)=\{x\in M: d(f^n(x),f^n(p))\xrightarrow[n\to+\infty]{}0\}
\end{align*}
Given $\beta>0$ we also consider the \emph{local stable set} of
size $\beta>0$, defined by
\begin{align*}
  W^s_\beta(p)=\{x\in M: d(f^n(x),f^n(p))\le\beta \text{ for all }n\ge0\}.
\end{align*}
The image of $W^s_\beta$ under the conjugacy $h$ is a neighborhood
of the origin inside $E^s$. It follows that the local stable set
is an embedded topological disk, with the same dimension as $E^s$.
Moreover, the orbits of the points in $W_\beta^s(p)$ actually
converges to the fixed point as time goes to infinity. Therefore,
$$
z \in W^s(p) \quad\Leftrightarrow \quad
 f^n(z) \in W^s_\beta(p) \text{ for some } n\ge 0.
$$

\subsection*{Stable manifold theorem}

The stable manifold theorem asserts that $W_\beta^s(p)$ is
actually a smooth embedded disk, with the same order of
differentiability as $f$ itself, and  it is tangent to $E^s$ at
the point $p$. It follows that $W^s(p)$ is a smooth submanifold,
injectively immersed in $M$. In general, $W^s(p)$ is not embedded
in $M$: in many cases it has self-accumulation points. For these
reasons one also refers to $W^s(p)$ and $W^s_\beta(p)$ as stable
\emph{manifolds} of $p$. Unstable manifolds are defined
analogously, replacing the transformation by its inverse.

\subsection*{Local stability}

We call \emph{index} of a diffeomorphism $f$ at a hyperbolic fixed
point $p$ the  index of the linear part, that is, the number of
eigenvalues of $Df_p$ with negative real part. By the
Hartman-Grobman theorem and previous comments on linear systems,
two diffeomorphisms are locally conjugate near hyperbolic fixed
points if and only if the stable indices and they both
preserve/reverse orientation. In other words, the index together
with the sign of the Jacobian determinant form a complete set of
invariants for local topological conjugacy.

Let $g$ be any diffeomorphism $C^1$-close to $f$. Then $g$ has a
unique fixed point $p_g$ close to $p$, and this fixed point is
still hyperbolic. Moreover, the stable indices and the
orientations of the two diffeomorphisms at the corresponding fixed
points coincide, and so they are locally conjugate. This is called
\emph{local stability} near of diffeomorphisms hyperbolic fixed
points. The same kind of result holds for flows near hyperbolic
equilibria.

\section{Hyperbolic behavior: examples}\label{sec:semi-global-hyperb}

Now let us review some key examples of (semi)global
hyperbolic dynamics. Thorough descriptions are available in
e.g. \cite{PM82}, \cite{KH95} and \cite{PT93}.

\subsection*{A linear torus automorphism}

Consider the linear transformation $A:\RR^2\to\RR^2$ given by the
following matrix, relative to the canonical base of the plane:
\begin{align*}
\begin{pmatrix}
  2 & 1 \\ 1 & 1
\end{pmatrix}.
\end{align*}
The $2$-dimensional torus $\TT^2$ is the quotient
$\RR^2/\ZZ^2$ of the plane by the equivalence relation
$$
(x_1,y_1) \sim (x_2,y_2) \quad\Leftrightarrow \quad
(x_1-x_2, y_1-y_2) \in \ZZ^2.
$$
Since $A$ preserves the lattice $\ZZ^2$ of integer vectors, that
is, since $A(\ZZ^2)=\ZZ^2$, the linear transformation defines an
invertible map $f_A:\TT^2\to\TT^2$ in the quotient space, which is
an example of linear automorphism of $\TT^2$. We call affine line
in $\TT^2$ the projection under the quotient map of any  affine
line in the plane.

The linear transformation $A$ is hyperbolic, with eigenvalues $0 <
\lambda_1 < 1 < \lambda_2$, and the corresponding eigenspaces
$E^1$ and $E^2$ have irrational slope. For each point $z\in
\TT^2$, let $W_i(z)$ denote the affine line through $z$ and having
the direction of $E^i$, for $i=1$, $2$:
\begin{itemize}
\item distances along $W_1(z)$ are multiplied by $\lambda_1<1$ under
      forward iteration of $f_A$
\item distances along $W_2(z)$ are multiplied by $1/\lambda_2<1$ under
      backward iteration of $f_A$.
\end{itemize}
Thus we call $W_1(z)$ \emph{stable manifold} and $W_2(z)$
\emph{unstable manifold} of $z$ (notice we are not assuming $z$ to
be periodic). Since the slopes are irrational, stable and unstable
manifolds are dense in the whole torus. From this fact one can
deduce that the periodic points of $f_A$ form a dense subset of
the torus, and that there exist points whose trajectories are
dense in $T^2$. The latter property is called \emph{transitivity}.

An important feature of this systems is that its behavior is
(globally) stable under small perturbations: given any
diffeomorphism $g:\TT^2\to\TT^2$ sufficiently $C^1$-close to
$f_A$, there exists a homeomorphism $h:\TT^2\to\TT^2$ such that $h
\circ g = f_A \circ h$. In particular, $g$ is also transitive and
its periodic points form a dense subset of $\TT^2$.

\subsection*{The Smale horseshoe}

\begin{figure}[htpb]
  \centering
   \psfrag{a}{$A$}\psfrag{b}{$B$}\psfrag{c}{$C$}\psfrag{d}{$D$}
   \psfrag{fa}{$f(A)$} \psfrag{fb}{$f(B)$}\psfrag{fc}{$f(C)$}
   \psfrag{0}{$0$}\psfrag{1}{$1$}
  \includegraphics[width=8cm]{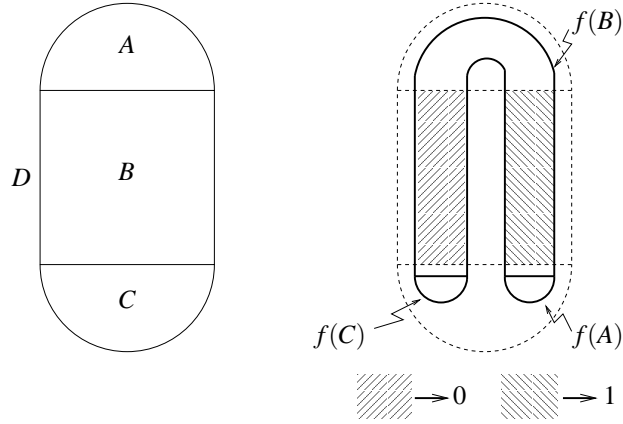}
  \caption{Horseshoe map\label{fig:horsesh-map}}
\end{figure}

Consider a stadium shaped region $D$ in the plane divided into
three subregions, as depicted in Figure~\ref{fig:horsesh-map}: two
half disks, $A$ and $C$, and a square, $B$. Next, consider a map
$f:D\to D$ mapping $D$ back inside itself as described in
Figure~\ref{fig:horsesh-map}: the intersection between $B$ and
$f(B)$ consists of two rectangles, $R_0$ and $R_1$, and $f$ is
affine on the pre-image of these rectangles, contracting the
horizontal direction and expanding the vertical direction.

The set $\Lambda=\cap_{n\in\ZZ}f^n(B)$, formed by all the points
whose orbits never leave the square $B$ is totally disconnected,
in fact, it is the product of two Cantor sets. A description of
the dynamics on $\Lambda$ may be obtained through the following
coding of orbits. For each point $z\in\Lambda$ and every time
$n\in\ZZ$ the iterate $f^n(z)$ must belong to either $R_0$ or
$R_1$. We call \emph{itinerary} of $z$ the sequence
$\{s_n\}_{n\in\ZZ}$ with values in the set $\{0, 1\}$ defined by
$f^n(z) \in R_{s_n}$ for all $n\in\ZZ$. The itinerary map
$$
\Lambda \to \{0,1\}^\ZZ, \quad z \mapsto \{s_n\}_{n\in\ZZ}
$$
is a homeomorphism, and conjugates $f$ restricted to $\Lambda$ to
the so-called \emph{shift map} defined on the space of sequences
by
$$
\{0,1\}^\ZZ \to \{0,1\}^\ZZ, \quad \{s_n\}_{n\in\ZZ}\mapsto \{s_{n+1}\}_{n\in\ZZ}.
$$
Since the shift map is transitive, and its periodic points form a
dense subset of the domain, it follows that the same is true for
the horseshoe map on $\Lambda$.

From the definition of $f$ we get that distances along horizontal
line segments through points of $\Lambda$ are contracted at a
uniform rate under forward iteration and, dually, distances along
vertical line segments through points of $\Lambda$ are contracted
at a uniform rate under backward iteration. Thus, horizontal line
segments are local stable sets and vertical line segments are
local unstable sets for the points of $\Lambda$.

A striking feature of this system is the stability of its
dynamics: given any diffeomorphism $g$ sufficiently $C^1$-close to
$f$, its restriction to the set $\Lambda_g=\cap_{n\in\ZZ}g^n(B)$
is conjugate to the restriction of $f$ to the set
$\Lambda=\Lambda_f$ (and, consequently, is conjugate to the shift
map). In addition, each point of $\Lambda_g$ has local stable and
unstable sets which are smooth curve segments, respectively,
approximately horizontal and approximately vertical.

\subsection*{The solenoid attractor}

\begin{figure}[ht]
  \centering
  \psfrag{a}{$\SS^1\times\DD$}
  \psfrag{b}{$f(\SS^1\times\DD)$}
  \psfrag{c}{$\{\theta\}\times\DD$}
  \includegraphics[width=6cm]{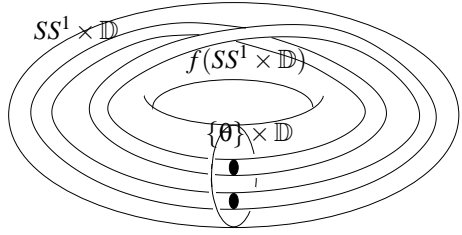}
  \caption{The solenoid attractor\label{fig:soleno-attract}}
  \end{figure}

The \emph{solid torus} is the product space $\SS^1\times\DD$,
where $\SS^1=\RR/\ZZ$ is the circle and $\DD=\{z\in\CC: |z|<1\}$
is the unit disk in the complex plane. Consider the map
$f:\SS^1\times\DD\to\SS^1\times\DD$ given by
\begin{align*}
  (\theta,z) \mapsto (2\theta,\alpha z + \beta e^{i\theta/2}),
\end{align*}
$\theta\in\RR/\ZZ$ and $\alpha$, $\beta\in\RR$ with
$\alpha+\beta<1$. The latter condition ensures that the image
$f(\SS^1\times\DD)$ is strictly contained in $\SS^1\times\DD$.
Geometrically, the image is a long thin domain going around the
solid torus twice, as described in
Figure~\ref{fig:soleno-attract}. Then, for any $n \ge 1$, the
corresponding iterate $f^n(\SS^1\times\DD)$ is an increasingly
thinner and longer domain that winds $2^k$ times around
$\SS^1\times\DD$. The maximal invariant set
$$
\Lambda=\cap_{n \ge 0} f^n(\SS^1\times\DD)
$$
is called \emph{solenoid attractor}. Notice that the forward orbit
under $f$ of every point in $\SS^1\times\DD$ accumulates on
$\Lambda$. One can also check that the restriction of $f$ to the
attractor is transitive, and the set of periodic points of $f$ is
dense in $\Lambda$.

In addition $\Lambda$ has a dense subset of periodic orbits
and also a dense orbit. Moreover every point in a
neighborhood of $\Lambda$ converges to $\Lambda$ and this is
why this set is called an \emph{attractor}.

\section{Hyperbolic sets}\label{sec:uniformly-hyperb-set}

The notion we are now going to introduce distillates the
crucial feature common to the examples presented previously.
A detailed presentation is given in e.g. \cite{PM82} and
\cite{Sh87}.  Let $f:M\to M$ be a diffeomorphism on a
manifold $M$.  A compact invariant set $\Lambda\subset M$ is
a \emph{hyperbolic set} for $f$ if the tangent bundle over
$\Lambda$ admits a decomposition
\begin{equation*}
T_\Lambda M = E^u \oplus E^s,
\end{equation*}
invariant under the derivative and such that $\|Df^{-1} \mid
E^u\|<\lambda$ and $\|Df \mid E^s \|< \lambda$ for some constant
$\lambda<1$ and some choice of a Riemannian metric on the
manifold. When it exists, such a decomposition is necessarily
unique and continuous. We call $E^s$ the stable bundle and $E^u$
the unstable bundle of $f$ on the set $\Lambda$.

The definition of hyperbolicity for an invariant set of a smooth
flow containing no equilibria is similar, except that one asks
for an invariant decomposition $T_\Lambda M = E^u \oplus E^0 \oplus E^s$,
where $E^u$ and $E^s$ are as before and $E^0$ is a line bundle
tangent to the flow lines. An invariant set that contains equilibria
is hyperbolic if and only it consists of a finite number of points,
all of them hyperbolic equilibria. 

\subsection*{Cone fields}

The definition of hyperbolic set is difficult to use in concrete 
situations, because, in most cases, one does not know the stable and
unstable bundles explicitly. Fortunately, to prove that an invariant
set is hyperbolic it suffices to have some approximate knowledge of
these invariant subbundles. That is the contents of the invariant cone
field criterion: a compact invariant set is hyperbolic if and only
if there exists some continuous (not necessarily invariant)
decomposition $T_\Lambda M= E^1\oplus E^2$ of the tangent bundle,
some constant $\lambda<1$, and some cone field around $E^1$
$$
 C^1_a(x)=\{ v=v_1+v_2\in E^1_x \oplus E^2_x: \|v_2\|\le a\|v_1\|\},
 \quad x\in\Lambda
$$
which is
\begin{itemize}
\item[(a)] forward invariant: $Df_x(C^1_a(x))\subset C^1_{\lambda a}(f(x))$ and
\item[(b)] expanded by forward iteration: $\|Df_x(v)\| \ge \lambda^{-1} \|v\|$
      for every $v\in C^1_a(x)$
\end{itemize}
and there exists a cone field $C^2_b(x)$ around $E^2$ which is
backward invariant and expanded by backward iteration.

\subsection*{Robustness}

An easy, yet very important consequence is that hyperbolic sets
are robust under small modifications of the dynamics. Indeed,
suppose $\Lambda$ is a hyperbolic set for $f:M\to M$, and let
$C_a^1(x)$ and $C_b^2(x)$ be invariant cone fields as above. The
(non-invariant) decomposition $E^1\oplus E^2$ extends continuously
to some small neighborhood $U$ of $\Lambda$, and then so do the
cone fields. By continuity, conditions (a) and (b) above remain
valid on $U$, possibly for a slightly larger constant $\lambda$.
Most important, they also remain valid when $f$ is replaced by any
other diffeomorphism $g$ which is sufficiently $C^1$-close to it.
Thus, using the cone field criterion once more, every compact set
$K\subset U$ which is invariant under $g$ is a hyperbolic set for
$g$.

\subsection*{Stable manifold theorem}

Let $\Lambda$ be a hyperbolic set for a diffeomorphism $f:M\to M$.
Assume $f$ is of class $C^k$. Then there exist $\epsilon_0>0$ and
$0<\lambda<1$ and, for each $0<\epsilon\le\epsilon_0$ and
$x\in\Lambda$, the \emph{local stable manifold of size $\epsilon$}
$$
W^s_\epsilon(x)=\{ y\in M: \dist(f^n(y),f^n(x))\le\epsilon
 \text{ for all } n \ge 0 \}
$$
and the  \emph{local unstable manifold of size $\epsilon$}
$$
W^u_\epsilon(x)=\{ y\in M: \dist(f^{-n}(y),f^{-n}(x))\le\epsilon
 \text{ for all } n \ge 0 \}
$$
are $C^k$ embedded disks, tangent at $x$ to $E^s_x$ and $E^u_x$,
respectively, and satisfying
\begin{itemize}
\item $f(W^s_\epsilon(x))\subset W^s_\epsilon(f(x))$ and
$f^{-1}(W^u_\epsilon(x))\subset W^u_\epsilon(f^{-1}(x))$;
\item $\dist(f(x),f(y))\le \lambda \dist(x,y)$ for all $y\in W^s_\epsilon(x)$
\item $\dist(f^{-1}(x),f^{-1}(y))\le \lambda \dist(x,y)$ for all $y\in
W^u_\epsilon(x)$
\item $W_\epsilon^s(x)$ and $W_\epsilon^u(x)$ vary continuously with the
point $x$, in the $C^k$ topology.
\end{itemize}
Then, the \emph{global stable and unstable manifolds} of $x$,
$$
  W^s(x)=\bigcup_{n\ge0} f^{-n}\big(W^s_\epsilon(f^n(x))\big)
  \quad\text{and}\quad
  W^u(x)=\bigcup_{n\ge0} f^{n}\big(W^u_\epsilon(f^{-n}(x))\big),
$$
are smoothly immersed submanifolds of $M$, and they are
characterized by
\begin{align*}
  W^s(x)&=\{y\in
  M:\dist(f^n(y),f^n(x))\to 0 \text{ as } n\to\infty\}
  \\
  W^u(x)&=\{y\in
  M:\dist(f^{-n}(y),f^{-n}(x)) \to 0 \text{ as } n\to\infty\}.
\end{align*}

\subsection*{Shadowing property}

This crucial property of hyperbolic sets means that possible small
``errors'' in the iteration of the map close to the set are, in
some sense, unimportant: to the resulting ``wrong'' trajectory,
there corresponds a nearby genuine orbit of the map. Let us give 
the formal statement. Recall that a hyperbolic set is compact,
by definition.

Given $\delta>0$, a $\delta$-\emph{pseudo-orbit} of $f:M\to M$ 
is a sequence $\{x_n\}_{n\in\ZZ}$ such that
$$
\dist(x_{n+1},f(x_n))\le\delta \quad\text{for all $n\in\ZZ$.}
$$
Given $\epsilon>0$, one says that a pseudo-orbit is
$\epsilon$-\emph{shadowed} by the orbit of a point $z\in M $ if
$\dist(f^n(z),x_n)\le\epsilon$ for all $n\in\ZZ$. The
\emph{shadowing lemma} says that for any $\epsilon>0$ one can find
$\delta>0$ and a neighborhood $U$ of the hyperbolic set $\Lambda$
such that every $\delta$-pseudo-orbit in $U$ is
$\epsilon$-shadowed by some orbit in $U$. Assuming $\epsilon$ is
sufficiently small, the shadowing orbit is actually unique.

\subsection*{Local product structure}

In general, these shadowing orbits need not be inside th
hyperbolic set $\Lambda$. However, that is indeed the case if
$\Lambda$ is a \emph{maximal invariant set}, that is, if it admits
some neighborhood $U$ such that $\Lambda$ coincides with the set
of points whose orbits never leave $U$:
$$
\Lambda=\bigcap_{n\in\ZZ} f^{-n}(U).
$$
A hyperbolic set is a maximal invariant set if and only if it has
the local product structure property stated in the next paragraph.

Let $\Lambda$ be a hyperbolic set and $\epsilon$ be small. If $x$
and $y$ are nearby points in $\Lambda$ then the local stable
manifold of $x$ intersects the local unstable manifold of $y$ at a
unique point, denoted $[x,y]$, and this intersection is
transverse. This is because the local stable manifold and the
local unstable manifold of every point are transverse, and these
local invariant manifolds vary continuously with the point. We say
that $\Lambda$ has \emph{local product structure} if there exists
$\delta>0$ such that $[x,y]$ belongs to $\Lambda$ for every $x$,
$y\in\Lambda$ with $\dist(x,y)<\delta$.

\subsection*{Stability}

The shadowing property may also be used to prove that hyperbolic
sets are stable under small perturbations of the dynamics:
if $\Lambda$ is a hyperbolic set for $f$ then for any $C^1$-close
diffeomorphism $g$ there exists a hyperbolic set $\Lambda_g$
close to $\Lambda$ and carrying the same dynamical behavior.

The key observation is that every orbit $f^n(x)$ of $f$ inside
$\Lambda$ is a $\delta$-pseudo-orbits for $g$ in a neighborhood $U$,
where $\delta$ is small if $g$ is close to $f$ and, hence, it is
shadowed by some orbit $g^n(z)$ of $g$. The correspondence $h(x)=z$
thus defined is injective and continuous.

For any diffeomorphism $g$ close enough to $f$, the orbits of $x$
in the maximal $g$-invariant set $\Lambda_g(U)$ inside $U$ are
pseudo-orbits for $f$. Therefore the shadowing property above
enables one to bijectively associate $g$-orbits of $\Lambda_g(U)$
to $f$-orbits in $\Lambda$. This provides a homeomorphism
$h:\Lambda_g(U)\to\Lambda$ which conjugates $g$ and $f$ on the
respective hyperbolic sets: $f\circ h= h \circ g$. Thus
\emph{hyperbolic maximal sets are structurally stable}: the
persistent dynamics in a neighborhood of these sets is the same
for all nearby maps.

If $\Lambda$ is a hyperbolic maximal invariant set for $f$ then
its hyperbolic continuation for any nearby diffeomorphism $g$ is
also a maximal invariant set for $g$.

\subsection*{Symbolic dynamics}

The dynamics of hyperbolic sets can be described through a symbolic
coding obtained from a convenient discretization of the phase space.
In a few words, one partitions the set into a finite number of 
subsets and assigns to a generic point in the hyperbolic set its
itinerary with respect to this partition. Dynamical properties can
then be read out from a shift map in the space of (admissible)
itineraries. The precise notion involved is that of Markov partition.

A set $R\subset\Lambda$ is a \emph{rectangle} if $[x,y]\in R$
for each $x,y\in R$. A rectangle is \emph{proper} if it is the
closure of its interior relative to $\Lambda$.
A \emph{Markov partition} of a hyperbolic set $\Lambda$ is a
cover $\cR=\{R_1,\dots,R_m\}$ of $\Lambda$ by proper rectangles
with pairwise disjoint interiors, relative to $\Lambda$, and such
$$
W^u(f(x))\cap R_j\subset f(W^u(x)\cap R_i)
 \quad\text{and}\quad
f(W^s(x)\cap R_i)\subset W^s(f(x))\cap R_j
$$
for every $x\in\inter_\Lambda(R_i)$ with $f(x)\in\inter_\Lambda
(R_j)$. The key fact is that \emph{any maximal hyperbolic set
$\Lambda$ admits Markov partitions with arbitrarily small diameter}.

Given a Markov partition $\cR$ with sufficiently small diameter,
and a sequence $\textbf{j}=(j_n)_{n\in\ZZ}$ in $\{1, \ldots, m\}$,
there exists at most one point $x=h(\textbf{j})$ such that
$$
f^n(x) \in R_{j_n} \quad\text{for each } n\in\ZZ.
$$
We say that $\textbf{j}$ is admissible if such a point $x$ does
exist and, in this case, we say $x$ admits $\textbf{j}$ as an
itinerary. It is clear that $f\circ h = h \circ \sigma$, where
$\sigma$ is the shift (left-translation) in the space of admissible
itineraries. The map $h$ is continuous and surjective, and it is 
injective on the residual set of points whose orbits never hit 
the boundaries (relative to $\Lambda$) of the Markov rectangles.

\section{Uniformly hyperbolic systems}\label{sec:uniformly-hyperb-sys}

A diffeomorphism $f:M\to M$ is \emph{uniformly hyperbolic},
or satisfies the \emph{Axiom A}, if the non-wandering set
$\Omega(f)$ is a hyperbolic set for $f$ and the set
$\Per(f)$ of periodic points is dense in $\Omega(f)$.  There
is an analogous definition for smooth flows $f^t:M\to M$,
$t\in\RR$. The reader can find the technical details in
e.g. \cite{KH95}, \cite{PM82} and \cite{Sh87}.

\subsection*{Dynamical decomposition}

The so-called ``spectral'' decomposition theorem of Smale
allows for the global dynamics of a hyperbolic
diffeomorphism to be decomposed into elementary building
blocks. It asserts that the non-wandering set splits into a
finite number of pairwise disjoint \emph{basic pieces} that
are compact, invariant, and dynamically indecomposable. More
precisely, the non-wandering set $\Omega(f)$ of a uniformly
hyperbolic diffeomorphism $f$ is a finite pairwise disjoint
union
$$
\Omega(f)=\Lambda_1 \cup \cdots \cup \Lambda_N
$$
of $f$-invariant, transitive sets $\Lambda_i$, that are compact and
maximal invariant sets. Moreover, the $\alpha$-limit set of every
orbit is contained in some $\Lambda_i$ and so is the $\omega$-limit set.

\subsection*{Geodesic flows on surfaces with negative curvature}

Historically, the first important example of uniform hyperbolicity 
was the geodesic flow $G^t$ on Riemannian manifolds of negative 
curvature $M$. This is defined as follows.

Let $M$ be a compact Riemannian manifold. Given any tangent vector
$v$, let $\gamma_v:\RR\to TM$ be the geodesic with initial condition
$v=\gamma_v(0)$. We denote by $\dot\gamma_v(t)$ the velocity vector
at time $t$. Since $\|\dot\gamma_v(t)\| = \|v\|$ for all $t$, it 
is no restriction to consider only unit vectors. 
There is an important volume form on the unit tangent bundle,
given by the product of the volume element on the manifold by the
volume element induced on each fiber by the Riemannian metric.
By integration of this form, one obtains the \emph{Liouville mesure}
on the unit tangent bundle, which is a finite measure if the 
manifold itself has finite volume (including the compact case).
The \emph{geodesic flow} is the flow $G^t:T^1M\to T^1M$ on the unit
tangent bundle $T^1M$ of the manifold, defined by
$$
G^t(v)=\dot\gamma_v(t).
$$
An important feature is that this flow leaves invariant the Liouville
measure. By Poincar\'e recurrence, this implies that $\Omega(G)=T^1 M$.

A major classical result in Dynamics, due to Anosov, states that
\emph{if $M$ has negative sectional curvature then this measure is
ergodic for the flow}. That is, any invariant set has zero or full
Liouville measure. The special case when $M$ is a surface, had been
dealt before by Hedlund and Hopf.

The key ingredient to this theorem is to prove that the geodesic 
flow is uniformly hyperbolic, in the sense we have just described,
when the sectional curvature is negative. In the surface case,
the stable and unstable invariant subbundles are differentiable,
which is no longer the case in general in higher dimensions.
This formidable obstacle was overcome by Anosov through showing
that the corresponding invariant foliations retain, nevertheless,
a weaker form of regularity property, that suffices for the proof.
Let us explain this.

\subsection*{Absolute continuity of foliations}

The invariant spaces $E^s_x$ and $E^u_x$ of a hyperbolic system
depend continuously, and even H\" older continuously, on the 
base point $x$. However, in general this dependence is not
differentiable, and this fact is at the origin of several important
difficulties. Related to this, the families of stable and unstable
manifolds are, usually, not differentiable foliations:
although the leaves themselves are as smooth as the dynamical system
itself, the holonomy maps often fail to be differentiable.
By holonomy maps we mean the projections along the leaves between
two given cross-sections to the foliation.    
 
However, Anosov and Sinai observed that if the system is at least 
twice differentiable then these foliations are \emph{absolutely
continuous}: their holonomy maps send zero Lebesgue measure sets 
of one cross-section to zero Lebesgue measure sets of the other 
cross-section. This property is crucial for proving that any smooth
measure which is invariant under a twice differentiable hyperbolic
system is ergodic. For dynamical systems that are only once 
differentiable the invariant foliations may fail to be absolutely
continuous. Ergodicity still is an open problem.

\subsection*{Structural stability}

A dynamical system is \emph{structurally stable} if it is equivalent
to any other system in a $C^1$ neighborhood, meaning that there
exists a global homeomorphism sending orbits of one to orbits of the
other and preserving the direction of time. More generally, replacing
$C^1$ by $C^r$ neighborhoods, any $r\ge 1$, one obtains the notion
of $C^r$ structural stability. Notice that, in principle, this 
property gets weaker as $r$ increases.

The Stability Conjecture of Palis-Smale proposed a complete geometric
characterization of this notion: for any $r\ge 1$, 
\emph{$C^r$ structurally stable systems should coincide with the
hyperbolic systems having the property of strong transversality},
that is, such that the stable and unstable manifolds of any points
in the non-wandering set are transversal. In particular, this
would imply that the property of $C^r$ structural stability does
not really depend on the value of $r$. 

That hyperbolicity and strong transversality suffice for structural
stability was proved in the 1970's by Robbin, de Melo, Robinson.
It is comparatively easy to prove that strong transversality is also
necessary. Thus, the heart of the conjecture is to prove that 
structurally stable systems must be hyperbolic. This was achieved
by Ma\~n\'e in the 1980's, for $C^1$ diffeomorphisms, and extended
about ten years later by Hayashi for $C^1$ flows. Thus
\emph{a $C^1$ diffeomorphism, or flow, on a compact manifold is
structurally stable if and only if it is uniformly hyperbolic and
satisfies the strong transversality condition.}

\subsection*{$\Omega$-stability}

A weaker property, called \emph{$\Omega$-stability} is defined
requiring equivalence only restricted to the non-wandering set.
The $\Omega$-Stability Conjecture of Palis-Smale claims that, for
any $r\ge 1$, \emph{$\Omega$-stable systems should coincide with
the hyperbolic systems with no cycles}, that is, such that no
basic pieces in the spectral decomposition are cyclically related
by intersections of the corresponding stable and unstable sets.

The $\Omega$-stability theorem of Smale states that these
properties are sufficient for $C^r$ $\Omega$-stability. Palis showed
that the no-cycles condition is also necessary. Much later, based 
on Ma\~n\'e's aforementioned result, he also proved that for $C^1$
diffeomorphisms hyperbolicity is necessary for $\Omega$-stability.
This was extended to $C^1$ flows by Hayashi in the 1990's.

\section{Attractors and physical measures}\label{sec:ergodic-theory-unifo}

A hyperbolic basic piece $\Lambda_i$ is a \emph{hyperbolic attractor}
if the stable set
$$
W^s(\Lambda_i)=\{x\in M : \omega(x)\subset \Lambda_i\}
$$
contains a neighborhood of $\Lambda_i$. In this case we call
$W^s(\Lambda_i)$ the \emph{basin} of the attractor
$\Lambda_i$\,, and denote it $B(\Lambda_i)$. When the
uniformly hyperbolic system is of class $C^2$, a basic piece
is an attractor if and only if its stable set has positive
Lebesgue measure. Thus, the union of the basins of all
attractors is a full Lebesgue measure subset of $M$. This
remains true for a residual (dense $G_\delta$) subset of
$C^1$ uniformly hyperbolic diffeomorphisms and flows.

The following fundamental result, due to Sinai, Ruelle,
Bowen shows that, no matter how complicated it may be, the
behavior of typical orbits in the basin of a hyperbolic
attractor is well-defined at the statistical level:
\emph{any hyperbolic attractor $\Lambda$ of a $C^2$ diffeomorphism
(or flow) supports a unique invariant probability measure $\mu$ such
that
\begin{equation}\label{eq.SinaiRuelleBowen}
 \lim_{n\to\infty} \frac{1}{n} \sum_{j=0}^{n-1} \varphi(f^j(z))
  = \int \varphi\,d\mu
\end{equation}
for every continuous function $\varphi$ and Lebesgue almost every
point $x\in B(\Lambda)$.
}
The standard reference here is~\cite{Bo75}.

Property \eqref{eq.SinaiRuelleBowen} also means that the
Sinai-Ruelle-Bowen measure $\mu$ may be ``observed'': the weights
of subsets may be found with any degree of precision, as the 
sojourn-time of any orbit picked ``at random'' in the basin
of attraction:
$$
\mu(V)=\text{ fraction of time the orbit of $z$ spends in $V$}
$$
for typical subsets $V$ of $M$ (the boundary of $V$ should have
zero $\mu$-measure), and for Lebesgue almost any point $z\in
B(\Lambda)$. For this reason $\mu$ is called a \emph{physical
measure}.

It also follows from the construction of these physical measures
on hyperbolic attractors that they depend continuously on the
diffeomorphism (or the flow). This \emph{statistical stability} is
another sense in which the asymptotic behavior is stable under
perturbations of the system, distinct from structural stability.
 
There is another sense in which this measure is ``physical''
and that is that $\mu$ is the zero-noise limit of the
stationary measures associated to the stochastic processes
obtained by adding small random noise to the system. The idea
is to replace genuine trajectories by ``random orbits'' $(z_n)_n$,
where each $z_{n+1}$ is chosen $\vep$-close to $f(z_n)$.
We speak of \emph{stochastic stability} if, for any continuous
function $\varphi$, the random time average
$$
 \lim_{n\to\infty} \frac{1}{n} \sum_{j=0}^{n-1} \varphi(z_j)
$$
is close to $\int\varphi\,d\mu$ for almost all choices of
the random orbit.

One way to construct such random orbits is through randomly perturbed
iterations, as follows. Consider a family of probability measures
$\nu_\vep$ in the space of diffeomorphisms, such that each $\nu_\vep$
is supported in the $\vep$-neighborhood of $f$. Then, for each initial
state $z_0$ define $z_{n+1}=f_{n+1}(z_n)$, where the diffeomorphisms
$f_n$ are independent random variables with distribution law $\nu_\vep$.
A probability measure $\eta_\vep$ on the basin $B(\Lambda)$ is
\emph{stationary}
if it satisfies
$$
\eta_\vep(E) = \int \eta_\vep(g^{-1}(E))\,d\nu_\vep(g).
$$
Stationary measures always exist, and they are often unique for
each small $\vep>0$. Then stochastic stability corresponds to
having $\eta_\vep$ converging weakly to $\mu$ when the noise level
$\vep$ goes to zero.

The notion of stochastic stability goes back to Kolmogorov and Sinai.
The first results, showing that uniformly hyperbolic systems are
stochastically stable, on the basin of each attractor, were proved in
the 1980's by Kifer and Young.

\begin{figure}[h]
\begin{center}
\psfrag{A}{$A$}
\psfrag{B}{$B$}
\psfrag{z}{$z$}
\includegraphics[height=1in]{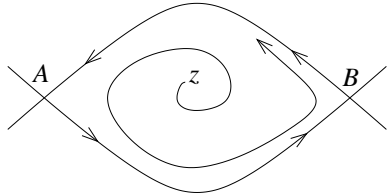}
\caption{\label{f.bowen} A planar flow with divergent time
averages}
\end{center}
\end{figure}

Let us point out that physical measures need not exist for
general systems.  A simple counter-example, attributed to Bowen, is
described in Figure~\ref{f.bowen}: time averages diverge over any of
the spiraling orbits in the region bounded by the saddle connections.
Notice that the saddle connections are easily broken by arbitrarily 
small perturbations of the flow. Indeed, no robust examples are known
of systems whose time-averages diverge on positive volume sets.

\section{Obstructions to hyperbolicity}\label{sec:obstruct-hyperb}

Although uniform hyperbolicity was originally intended to encompass a
residual or, at least, dense subset of all dynamical systems, it was
soon realized that this is not the case: many important examples fall
outside its realm. There are two main mechanisms
that yield robustly non-hyperbolic behavior, that is, whole open sets
of non-hyperbolic systems.

\subsection*{Heterodimensional cycles}

\begin{figure}[pht]
\begin{center}
 \psfrag{q}{$q$}
 \psfrag{p1}{$p_1$}
 \psfrag{p2}{$p_2$}
\includegraphics[height=1.5in]{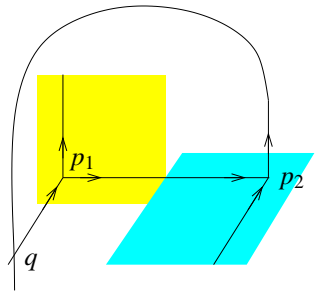}
\caption{\label{f.heterodimensional}A heterodimensional cycle}
\end{center}
\end{figure}

Historically, the first such mechanism was the coexistence of periodic
points with different Morse indices (dimensions of the unstable manifolds)
inside the same transitive set.
See Figure~\ref{f.heterodimensional}. This is how the first examples
of $C^1$-open subsets of non-hyperbolic diffeomorphisms were
obtained by Abraham, Smale on manifolds of dimension $d\ge 3$.
It was also the key in the constructions by Shub and Ma\~n\'e of
non-hyperbolic, yet robustly transitive diffeomorphisms, that is,
such that every diffeomorphism in a $C^1$ neighborhood has dense orbits.

For flows, this mechanism may assume a novel form, because
of the interplay between regular orbits and singularities
(equilibrium points).  That is, robust non-hyperbolicity may
stem from the coexistence of regular and singular orbits in
the same transitive set. The first, and very striking
example was the geometric Lorenz attractor proposed by
Afraimovich, Bykov, Shil'nikov and Guckenheimer, Williams to
model the behavior of the Lorenz equations, that we shall
discuss later.

\subsection*{Homoclinic tangencies}

\begin{figure}[pht]
\begin{center}
\psfrag{H}{$H$} \psfrag{p}{$p$}
\includegraphics[height=1.5in]{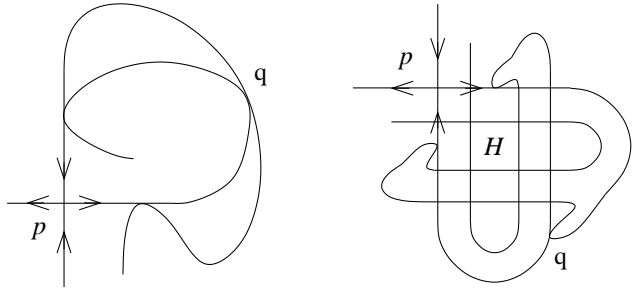}
\caption{\label{f.homoclinic} Homoclinic tangencies}
\end{center}
\end{figure}

Of course, heterodimensional cycles may exist only in dimension $3$ or
higher. The first robust examples of non-hyperbolic diffeomorphisms on
surfaces were constructed by Newhouse, exploiting the second of these
two mechanisms: homoclinic tangencies, or non-transverse intersections
between the stable and the unstable manifold of the same periodic point.
See Figure~\ref{f.homoclinic}.

It is important to observe that individual homoclinic
tangencies are easily destroyed by small perturbations of
the invariant manifolds.  To construct open examples of
surface diffeomorphisms with \emph{some} tangency, Newhouse
started from systems where the tangency is associated to a
periodic point inside an invariant hyperbolic set with rich
geometric structure. This is illustrated on the right hand
side of Figure~\ref{f.homoclinic}. His argument requires a
very delicate control of distortion, as well as of the
dependence of the fractal dimension on the
dynamics. Actually, for this reason, his construction is
restricted to the $C^r$ topology for $r\ge 2$. A very
striking consequence of this construction is that these open
sets exhibit \emph{coexistence of infinitely many periodic
  attractors}, for each diffeomorphism on a residual
subset. A detailed presentation of his result and
consequences is given in \cite{PT93}.

Newhouse's conclusions have been extended in two
ways. First, by Palis, Viana, for diffeomorphisms in any
dimension, still in the $C^r$ topology with $r\ge 2$. Then,
by Bonatti, D\'\i az, for $C^1$ diffeomorphisms in any
dimension larger or equal than $3$. The case of $C^1$
diffeomorphisms on surfaces remains open. As a matter of
fact, in this setting it is still unknown whether uniform
hyperbolicity is dense in the space of all diffeomorphisms.

\section{Partial hyperbolicity}\label{sec:robustn-weak-hyperb}

Several extensions of the theory of uniform hyperbolicity have been
proposed, allowing for more flexibility, while keeping the core idea:
splitting of the tangent bundle into invariant subbundles.
We are going to discuss more closely two such extensions. 

On the one hand, one may allow for one or more invariant
subbundles along which the derivative exhibits mixed
contracting/neutral/expanding behavior. This is generically
referred to as \emph{partial hyperbolicity}, and a standard
reference is the book \cite{HPS77}.  On the other
hand, while requiring all invariant subbundles to be either
expanding or contraction, one may relax the requirement of
uniform rates of expansion and contraction. This is usually
called \emph{non-uniform hyperbolicity}. A detailed
presentation of the fundamental results about this notion is
available e.g. in \cite{KH95}. In this section we
discuss the first type of condition. The second one will be
dealt with later.

\subsection*{Dominated splittings}

Let $f:M\to M$ be a diffeomorphism on a closed manifold $M$ and
$K$ be any $f$-invariant set. A continuous splitting $T_x M=
E_1(x)\oplus\cdots\oplus E_k(x)$, $x\in K$ of the tangent bundle
over $K$ is \emph{dominated} if it is invariant under the
derivative $Df$ and there exists $\ell\in\NN$ such that for every
$i<j$, every $x\in K$, and every pair of unit vectors $u\in
E_i(x)$ and $v\in E_j(x)$, one has
\begin{equation}\label{eq.domination}
\frac{\|Df^\ell_x\cdot u\|}{\|Df^\ell_x\cdot v\|}<\frac12
\end{equation}
and the dimension of $E_i(x)$ is independent of $x\in K$
for every $i\in\{ 1,\dots, k\}$. 
This definition may be formulated, equivalently, as follows: there
exist $C>0$ and $\lambda<1$ such that for
every pair of unit vectors $u\in E_i(x)$ and $v\in E_j(x)$,
one has
\begin{equation*}
\frac{\|Df^n_x\cdot u\|}{\|Df^n_x\cdot v\|}< C \lambda^n \quad\text{for all
$n\ge 1$}.
\end{equation*}

Let $f$ be a diffeomorphism and $K$ be an $f$-invariant set having
a dominated splitting $T_K M = E_1\oplus\cdots\oplus E_k$. We say
that the splitting and the set $K$ are
\begin{itemize}
\item\emph{partially hyperbolic} the derivative either contracts
uniformly $E_1$ or expands uniformly $E_k$: there exists $\ell\in\NN$
such that 
$$
\text{either } \|Df^\ell\mid E_1\| < \frac 12 \quad\text{or}\quad
\|(Df^\ell\mid E_k)^{-1}\| < \frac 12.
$$
 
\item\emph{volume hyperbolic} if the derivative either contracts
volume uniformly along $E_1$ or expands volume uniformly along $E_k$:
there exists $\ell\in\NN$ such that 
$$
\text{either } |\det(Df^\ell\mid E_1)| < \frac 12 \quad\text{or}\quad
|\det(Df^\ell\mid E_k)| > 2.
$$
\end{itemize}

The diffeomorphism $f$ is \emph{partially hyperbolic/volume hyperbolic}
if the ambient space $M$ is a partially hyperbolic/volume hyperbolic set
for $f$.

\subsection*{Invariant foliations}

An crucial geometric feature of partially hyperbolic systems is the
existence of invariant foliations tangent to uniformly expanding or 
uniformly contracting invariant subbundles: \emph{assuming the
derivative contracts $E^1$ uniformly, there exists a unique family
$\cF^s=\{\cF^s(x):x\in K\}$ of injectively $C^r$ immersed submanifolds
tangent to $E^1$ at every point of $K$, satisfying 
$f(\cF^s(x))=\cF^s(f(x))$ for all $x\in K$, and which are uniformly
contracted by forward iterates of $f$.}  This is called 
\emph{strong-stable foliation} of the diffeomorphism on $K$.
Strong-unstable foliations are defined in the same way, tangent to
the invariant subbundle $E_k$, when it is uniformly expanding.

As in the purely hyperbolic setting, a crucial ingredient in the ergodic
theory of partially hyperbolic systems is the fact that strong-stable
and strong-unstable foliations are absolutely continuous, if the 
system is at least twice differentiable.

\subsection*{Robustness and partial hyperbolicity}

Partially hyperbolic systems have been studied since the 1970's, most
notably by Brin, Pesin and Hirsch, Pugh, Shub. Over the last decade
they attracted much attention as the key to characterizing robustness
of the dynamics. More precisely, let $\Lambda$ be a maximal invariant
set of some diffeomorphism $f$:
$$
\Lambda = \bigcap_{n\in\ZZ} f^n(U)
 \quad\text{for some neighborhood $U$ of $\Lambda$.}
$$ 
The set $\Lambda$ is \emph{robust}, or \emph{robustly transitive}, if
its continuation $\Lambda_g = \cap_{n\in\ZZ} g^n(U)$ is transitive for
all $g$ in a neighborhood of $f$. There is a corresponding notion for
flows.

As we have already seen, hyperbolic basic pieces are
robust. In the 1970's, Ma\~n\'e observed that the converse
is also true when $M$ is a surface, but not anymore if the
dimension of $M$ is at least $3$.  Counter-examples in
dimension $4$ had been given before by Shub.  A series of
results of Bonatti, D\'\i az, Pujals, Ures in the 1990's
clarified the situation in all dimensions: robust sets
always admit some dominated splitting which is volume
hyperbolic; in general, this splitting needs not be
partially hyperbolic, except when the ambient manifold has
dimension $3$.

\subsection*{Lorenz-like strange attractors}

Parallel results hold for flows on $3$-dimensional manifolds.
The main motivation are the so-called Lorenz-like strange attractors,
inspired by the famous differential equations
\begin{equation}\label{eq.lorenz}
\begin{array}{ll}
\dot x  = - \sigma x + \sigma y & \qquad \sigma = 10 \\
\dot y  = \rho x - y - xz       & \qquad \rho = 28      \\
\dot z  = xy - \beta z          & \qquad \beta =8/3
\end{array}
\end{equation}
introduced by E. N. Lorenz in the early 1960's. 
Numerical analysis of these equations led Lorenz to realize that sensitive
dependence of trajectories on the initial conditions is ubiquitous among
dynamical systems, even those with simple evolution laws.

The dynamical behavior of \eqref{eq.lorenz} was first interpreted by means
of certain geometric models, proposed by Guckenheimer, Williams and
Afraimovich, Bykov, Shil'nikov in the 1970's, where the presence of strange
attractors, both sensitive and fractal, could be proved rigorously.
It was much harder to prove that the original equations~\eqref{eq.lorenz}
themselves have such an attractor. This was achieved just a few years ago,
by Tucker, by means of a computer assisted rigorous argument. 

An important point is that Lorenz-like attractors cannot be
hyperbolic, because they contain an equilibrium point
accumulated by regular orbits inside the attractor. Yet,
these strange attractors are robust, in the sense we defined
above.  A mathematical theory of robustness for flows in
$3$-dimensional spaces was recently developed by Morales,
Pacifico, and Pujals. In particular, this theory shows that
uniformly hyperbolic attractors and Lorenz-like attractors
are the only ones which are robust. Indeed, they prove that
\emph{any robust invariant set of a flow in dimension $3$ is
  singular hyperbolic}. Moreover, \emph{if the robust set
  contains equilibrium points then it must be either an
  attractor or a repeller}. A detailed presentation of this
and related results is given in \cite{AraPac07}.

An invariant set $\Lambda$ of a flow in dimension $3$ is
\emph{singular hyperbolic} if it is a partially hyperbolic
set with splitting $E^1\oplus E^2$ such that the derivative
is volume contracting along $E^1$ and volume expanding along
$E^2$.  Notice that one of the subbundles $E^1$ or $E^2$
must be one-dimensional, and then the derivative is,
actually, either norm contracting or norm expanding along
this subbundle.  Singular hyperbolic sets without equilibria
are uniformly hyperbolic: the $2$-dimensional invariant
subbundle splits as the sum of the flow direction with a
uniformly expanding or contracting one-dimensional invariant
subbundle.

\section{Non-uniform hyperbolicity - Linear theory}

In its linear form, the theory of non-uniform hyperbolicity
goes back to Lyapunov, and is founded on the multiplicative
ergodic theorem of Oseledets. Let us introduce the main
ideas, whose thorough development can be found in
e.g. \cite{CFS82}, \cite{KH95} and \cite{Man87}.

The \emph{Lyapunov exponents} of a sequence $\{A^n, n\ge
1\}$ of square matrices of dimension $d\ge 1$, are the values of
\begin{equation}\label{exponent1}
\lambda(v) = \limsup_{n\to\infty} \frac 1n \log\|A^n\cdot v\|
\end{equation}
over all non-zero vectors $v\in\RR^d$. For completeness, set
$\lambda(0)=-\infty$. It is easy to see that
$\lambda(cv)=\lambda(v)$ and
$\lambda(v+v')\le\max\{\lambda(v),\lambda(v')\}$ for any non-zero
scalar $c$ and any vectors $v$, $v'$. It follows that, given any
constant $a$, the set of vectors satisfying $\lambda(v)\le a$ is a
vector subspace. Consequently, there are at most $d$ Lyapunov
exponents, henceforth denoted by $\lambda_1<\cdots< \lambda_{k-1} <
\lambda_k$, and there exists a filtration $F_0\subset F_1
\subset \cdots \subset  F_{k-1}
\subset F_k=\RR^d$ into vector subspaces, such that
$$
\lambda(v)=\lambda_i \text{ for all } v\in F_i\setminus F_{i-1}
$$
and every $i=1, \ldots, k$ (write $F_0=\{0\}$). In particular, the
largest exponent is given by
\begin{equation}\label{exponent2}
\lambda_k = \limsup_{n\to\infty} \frac 1n \log\|A^n\|\,.
\end{equation}
One calls $\dim F_i-\dim F_{i-1}$ the \emph{multiplicity} of each
Lyapunov exponent $\lambda_i$.

There are corresponding notions for continuous families of
matrices $A^t$, $t\in (0,\infty)$, taking the limit as $t$ goes to
infinity in the relations \eqref{exponent1} and \eqref{exponent2}.

\subsection*{Lyapunov stability}

Consider the linear differential equation
\begin{equation}\label{stability1}
\dot v(t) = B(t)\cdot v(t)
\end{equation}
where $B(t)$ is a bounded function with values in the space of
$d\times d$ matrices, defined for all $t\in\RR$. The theory of
differential equations ensures that there exists a
\emph{fundamental matrix} $A^t$, $t\in\RR$ such that
$$
v(t)=A^t\cdot v_0
$$
is the unique solution of \eqref{stability1} with initial
condition $v(0)=v_0$.

If the Lyapunov exponents of the family $A^t$, $t>0$ are all
negative then the trivial solution $v(t)\equiv 0$ is
asymptotically stable, and even exponentially stable. The
stability theorem of A. M. Lyapunov asserts that, under an
additional regularity condition, stability is still valid for
non-linear perturbations
\begin{equation*}
w(t)=B(t)\cdot w + F(t,w)
\end{equation*}
with $\|F(t,w)\|\le \const \|w\|^{1+c}$, $c>0$. That is, the
trivial solution $w(t)\equiv 0$ is still exponentially
asymptotically stable.

The regularity condition means, essentially, that the limit in
\eqref{exponent1} does exist, even if one replaces vectors $v$ by
elements $v_1 \wedge \cdots \wedge v_l$ of any $l$th exterior
power of $\RR^d$, $1 \le l \le d$. By definition, the norm of an
$l$-vector $v_1 \wedge \cdots \wedge v_l$ is the volume of the
parallelepiped determined by the vectors $v_1$, \ldots, $v_k$.
This condition is usually tricky to check in specific situations.
However, the multiplicative ergodic theorem of V. I. Oseledets
asserts that, for very general matrix-valued stationary random
processes, regularity is an almost sure property.

\subsection*{Multiplicative ergodic theorem}

Let $f:M\to M$ be a measurable transformation, preserving some
measure $\mu$, and let $A:M\to \operatorname{GL}(d,\RR)$ be any
measurable function such that $\log\|A(x)\|$ is $\mu$-integrable.
The Oseledets theorem states that Lyapunov exponents exist for
the sequence $A^n(x)=A(f^{n-1}(x)) \cdots A(f(x)) \, A(x)$ for
$\mu$-almost every $x\in M$. More precisely, for $\mu$-almost
every $x\in M$ there exists $k=k(x)\in\{1, \ldots, d\}$, 
a filtration
$$
F^0_x\subset F^1_x \subset \cdots \subset F_x^{k-1} \subset F_x^k=\RR^d,
$$
and numbers $\lambda_1(x)<\cdots<\lambda_k(x)$ such that
\begin{equation*}
\lim_{n\to\infty} \frac 1n \log\|A^n(x)\cdot v\| = \lambda_i(x)
\end{equation*}
for all $v\in F_x^i\setminus F_x^{i-1}$ and $i\in\{1,\ldots, k\}$.
More generally, this conclusion holds for any vector bundle 
automorphism $\cV\to\cV$ over the transformation $f$, with
$A_x:\cV_x\to\cV_{f(x)}$ denoting the action of the automorphism
on the fiber of $x$.

The Lyapunov exponents $\lambda_i(x)$, and their number $k(x)$,
are measurable functions of $x$ and they are constant on orbits of
the transformation $f$. In particular, if the measure $\mu$ is
ergodic then $k$ and the $\lambda_i$ are constant on a full
$\mu$-measure set of points. The subspaces $F_x^i$ also depend
measurably on the point $x$ and are invariant under the
automorphism:
$$
A(x) \cdot F_x^i = F_{f(x)}^i.
$$
It is in the nature of things that, usually, these objects are
\emph{not} defined everywhere and they depend discontinuously
on the base point $x$.

When the transformation $f$ is invertible one obtains a stronger
conclusion, by applying the previous result also to the inverse
automorphism: assuming that $\log\|A(x)^{-1}\|$ is also
in $L^1(\mu)$, one gets that there exists a decomposition
$$
\cV_x=E_x^1\oplus\cdots\oplus E_x^k,
$$
defined at almost every point and such that $A(x) \cdot E_x^i =
E_{f(x)}^i$ and
\begin{equation*}
\lim_{n\to\pm\infty} \frac 1n \log\|A^n(x)\cdot v\| = \lambda_i(x)
\end{equation*}
for all $v\in E_x^i$ different from zero and all $i\in\{1, \ldots,
k\}$. These \emph{Oseledets subspaces} $E_x^i$ are related to the
subspaces $F_x^i$ through
$$
F_x^j = \oplus_{i=1}^j E_x^i.
$$
Hence, $\dim E_x^i = \dim F_x^i-\dim F_x^{i-1}$ is the
multiplicity of the Lyapunov exponent $\lambda_i(x)$.

The angles between any two Oseledets subspaces decay
sub-exponentially along orbits of $f$:
$$
\lim_{n\to\pm\infty} \frac 1n \log
\ang{\bigoplus_{i\in I} E_{f^n(x)}^i, \bigoplus_{j\notin I} E_{f^n(x)}^j} = 0
$$
for any $I \subset \{1, \ldots, k\}$ and almost every point.
These facts imply the regularity condition mentioned previously and,
in particular,
\begin{equation*}
\lim_{n\to\pm\infty} \frac 1n \log | \det A^n(x)| =
\sum_{i=1}^k \lambda_i(x) \dim E_x^i
\end{equation*}
Consequently, if $\det A(x)=1$ at every point then the sum of all Lyapunov
exponents, counted with multiplicity, is identically zero.

\section{Non-uniformly hyperbolic systems}\label{sec:non-uniformly-hyperb}

The Oseledets theorem applies, in particular, when $f:M\to M$ is a $C^1$
diffeomorphism on some compact manifold and $A(x)=Df_x$. Notice that
the integrability conditions are automatically satisfied, for any 
$f$-invariant probability measure $\mu$, since the derivative of $f$ and
its inverse are bounded in norm.

Lyapunov exponents yield deep geometric information on the dynamics of
the diffeomorphism, especially when they do not vanish. 
We call $\mu$ a \emph{hyperbolic measure} if all Lyapunov exponents are
non-zero at $\mu$-almost every point. By \emph{non-uniformly hyperbolic
system} we shall mean a diffeomorphism $f:M\to M$ together with some
invariant hyperbolic measure.

A theory initiated by Pesin provides fundamental geometric
information on this class of systems, especially
existence of stable and unstable manifolds at almost every
point which form absolutely continuous invariant
laminations. For most results, one needs the derivative
$Df$ to be H\"older continuous: there exists $c>0$ such that
$$
\|Df_x- Df_y\|\le\const \cdot d(x,y)^c.
$$
These notions extend to the context of flows essentially
without change, except that one disregards the invariant
line bundle given by the flow direction (whose Lyapunov
exponent is always zero). A detailed presentation can be
found in e.g. \cite{KH95}.

\subsection*{Stable manifolds}

An essential tool is the existence of invariant families of
local stable sets and local unstable sets, defined at
$\mu$-almost every point.  Assume $\mu$ is a hyperbolic
measure.  Let $E_x^u$ and $E_x^s$ be the sums of all
Oseledets subspaces corresponding to positive, respectively
negative, Lyapunov exponents, and let $\tau_x>0$ be a lower
bound for the norm of every Lyapunov exponent at $x$.

Pesin's stable manifold theorem states that, \emph{for $\mu$-almost every
$x\in M$, there exists a $C^1$ embedded disk $W_{loc}^s(x)$ tangent to
$E^s_x$ at $x$ and there exists $C_x>0$ such that
$$
\dist(f^n(y),f^n(x))\le C_x e^{-n\tau_x}\cdot\dist(y,x)
 \quad\text{for all $y\in W_{loc}^s(x)$.}
$$
}
Moreover, the family $\{W^s_{loc}(x)\}$ is invariant, in the sense that
$f(W^s_{loc}(x))\subset W^s_{loc}(f(x))$ for $\mu$-almost every $x$.
Thus, one may define global stable manifolds
$$
W^s(x)=\bigcup_{n=0}^\infty f^{-n}\big(W^s_{loc}(x)\big)
 \quad\text{for $\mu$-almost every $x$.}
$$
In general, the local stable disks $W^s(x)$ depend only measurably on $x$.
Another key difference with respect to the uniformly hyperbolic setting
is that the numbers $C_x$ and $\tau_x$ can not be taken independent 
of the point, in general. Likewise, one defines local and global unstable
manifolds, tangent to $E^u_x$ at almost every point.
Most important for the applications, both foliations, stable and unstable,
are absolutely continuous.

In the remaining sections we briefly present three major results in the
theory of non-uniform hyperbolicity: the entropy formula, abundance of
periodic orbits, and exact dimensionality of hyperbolic measures.

\subsection*{The entropy formula}

The entropy of a partition $\cP$ of $M$ is defined by
$$
h_\mu(f,\cP)=\lim_{n\to\infty} \frac 1n H_\mu(\cP^n),
$$
where $\cP^n$ is the partition into sets of the form $P=P_0\cap
f^{-1}(P_1) \cap \cdots \cap f^{-n}(P_n)$ with $P_j\in\cP$ and
$$
H_\mu(\cP^n) = \sum_{P\in\cP^n} -\mu(P) \log\mu(P).
$$
The \emph{Kolmogorov-Sinai entropy} $h_\mu(f)$ of the system is
the supremum of $h_\mu(f,\cP)$ over all partitions $\cP$ with
finite entropy. The Ruelle-Margulis inequality says that
$h_\mu(f)$ is bounded above by the averaged sum of the positive
Lyapunov exponents. A major result of the theorem, Pesin's entropy
formula, asserts that if the invariant measure $\mu$ is smooth
(for instance, a volume element) then the entropy actually coincides
with the averaged sum of the positive Lyapunov exponents
$$
h_\mu(f) = \int \big(\sum_{j=1}^k \max\{0,\lambda_j\}\big) d\mu.
$$
A complete characterization of the invariant measures for which
the entropy formula is true was given by F. Ledrappier and L. S.
Young.

\subsection*{Periodic orbits and entropy}

It was proved by A. Katok that periodic motions are always dense in the
support of any hyperbolic measure. More than that, assuming the measure
is non-atomic, there exist Smale horseshoes $H_n$ with topological 
entropy arbitrarily close to the entropy $h_\mu(f)$ of the system.
In this context, the \emph{topological entropy} $h(f,H_n)$ may be defined
as the exponential rate of growth
$$
\lim_{k\to\infty} \frac 1k \log\#\{x\in H_n: f^k(x)=x\}.
$$
of the number of periodic points on $H_n$.

\subsection*{Dimension of hyperbolic measures}

Another remarkable feature of hyperbolic measures is that they are
\emph{exact dimensional}: the pointwise dimension
$$
d(x) = \lim_{r\to 0} \frac{\log\mu(B_r(x))}{\log r}
$$
exists at almost every point, where $B_r(x)$ is the
neighborhood of radius $r$ around $x$.  This fact was proved
by L. Barreira, Ya.  Pesin, and J. Schmeling.  Note that
this means that the measure $\mu(B_r(x))$ of neighborhoods
scales as $r^{d(x)}$ when the radius $r$ is small.

\section{Future directions}

The theory of uniform hyperbolicity showed that dynamical
systems with very complex behavior may be amenable to a very
precise description of their evolution, especially in
probabilistic terms.  It was most successful in
characterizing structural stability, and also established a
paradigm of how general "chaotic" systems might be
approached.  A vast research program has been going on in
the last couple of decades or so, to try and build such a
global theory of complex dynamical evolution, where notions
such as partial and non-uniform hyperbolicity play a central
part.  The reader is referred to the bibliography,
especially the book \cite{Beyond}
for a review of much recent progress.


\bibliographystyle{plain}

\begin{thebibliography}{10}

\bibitem{AraPac07}
V.~Araujo and M.~J. Pacifico.
\newblock {\em Three Dimensional Flows}.
\newblock XXV Brazillian Mathematical Colloquium. IMPA, Rio de Janeiro, 2007.

\bibitem{Beyond}
C.~Bonatti, L.~J. D{\'{\i}}az, and M.~Viana.
\newblock {\em Dynamics beyond uniform hyperbolicity}, volume 102 of {\em
  Encyclopaedia of Mathematical Sciences}.
\newblock Springer-Verlag, 2005.

\bibitem{Bo75}
R.~Bowen.
\newblock {\em Equilibrium states and the ergodic theory of {A}nosov
  diffeomorphisms}, volume 470 of {\em Lect. Notes in Math.}
\newblock Springer Verlag, 1975.


\bibitem{CFS82}
I.~P. Cornfeld, S.~V. Fomin, and Ya.~G. Sina{\u\i}.
\newblock {\em Ergodic theory}, volume 245 of {\em Grundlehren der
  Mathematischen Wissenschaften [Fundamental Principles of Mathematical
  Sciences]}.
\newblock Springer-Verlag, 1982.

\bibitem{HPS77}
M.~Hirsch, C.~Pugh, and M.~Shub.
\newblock {\em Invariant manifolds}, volume 583 of {\em Lect. Notes in Math.}
\newblock Springer Verlag, 1977.

\bibitem{KH95}
A.~Katok and B.~Hasselblatt.
\newblock {\em Introduction to the modern theory of dynamical systems}.
\newblock Cambridge University Press, 1995.

\bibitem{Man87}
R.~Ma{\~{n}}{\'{e}}.
\newblock {\em Ergodic theory and differentiable dynamics}.
\newblock Springer Verlag, 1987.


\bibitem{PM82}
J.~Palis and W.~de~Melo.
\newblock {\em Geometric theory of dynamical systems. An introduction}.
\newblock Springer Verlag, 1982.

\bibitem{PT93}
J.~Palis and F.~Takens.
\newblock {\em Hyperbolicity and sensitive-chaotic dynamics at homoclinic
  bifurcations}.
\newblock Cambridge University Press, 1993.

\bibitem{Sh87}
M.~Shub.
\newblock {\em Global stability of dynamical systems}.
\newblock Springer Verlag, 1987.

\end{thebibliography}

\end{document}